\documentstyle[amssymb]{article}
\newtheorem{theorem}{Theorem}[section]
\newtheorem{lemma}[theorem]{Lemma}

\newtheorem{definition}[theorem]{Definition}

\newcommand{\be}{\begin{equation}}
\newcommand{\ee}{\end{equation}}
\newcommand{\beq}{\begin{equation}}
\newcommand{\enq}{\end{equation}}

\newcommand{\Rr}{{\Bbb{R}}}
\newcommand{\Nn}{{\Bbb{N}}}
\newcommand{\Cc}{{\Bbb{C}}}

\newcommand{\lek}{ \Lambda_{\eta,K}}

  \title{Spectral inclusion and spectral exactness for singular non-selfadjoint
         Sturm-Liouville problems}
\author{ B.M. Brown   \\
Department of Computer Science, \\University of Wales, Cardiff, 
PO Box 916, Cardiff CF2 3XF,  U.K.  
\and M. Marletta \\ Department of Mathematics and Computer Science, \\ 
University of Leicester, \\
 University Road, Leicester LE1 7RH  U.K.}
\begin{document}
\maketitle

\begin{abstract}
We consider the effect of regularization by interval truncation on the spectrum 
of a singular non-selfadjoint Sturm-Liouville operator. We present results on
spectral inclusion and spectral exactness for the cases where the singularity
is in Sims Case II or Sims Case III. For Sims Case I we present a test for
spectral inexactness, which can be used to detect when the interval truncation
process is generating spurious eigenvalues. Numerical results illustrate the
effectiveness of this test.
\end{abstract}

\noindent {\sc AMS(MOS) Subject classifications:} 34L, 34B20, 47A74, 65L15

 
\section{Introduction}\label{section:1} 
Over the last 30 years there has been considerable interest in numerical
solution of singular Sturm-Liouville problems, and in particular in the
development of automatic software for such problems: see, e.g., Bailey,
Gordon and Shampine \cite{kn:BSG}, Bailey, Garbow, Kaper and Zettl \cite{kn:BGZ},
Fulton and Pruess \cite{kn:FP} and Marletta and Pryce \cite{kn:MP}. The
software described in these papers usually uses an interval truncation
procedure to regularize problems posed either on infinite intervals, or
on finite intervals with singular behaviour of the coefficients near at
least one of the endpoints. Rigorous mathematical justification of the
validity of the interval truncation process, however, did not appear
(except for special cases) until the paper of Bailey, Everitt, Weidmann 
and Zettl \cite{kn:BEWZ} in 1993, which uses fundamental ideas from
Reed and Simon \cite{kn:RS} to develop conditions under which the spectra
of a sequence of regularized problems can 
(a) provide approximations to the whole spectrum of the original
 singular eigenvalue problem ({\em spectral inclusion});
(b) {\em not} yield approximations to any points which are not
 in the spectrum of the original singular problem ({\em spectral exactness}).
All of this work is for selfadjoint problems only.

Non-selfadjoint singular problems are also very important. They arise when
the complex scaling method is used to find resonances of a selfadjoint
problem (for a review see \cite{kn:Hislop}) and also, more classically,
in the study of hydrodynamic stability, where the spectra of the Orr-Sommerfeld
and related equations are often studied over infinite intervals. Recent applications 
to the study of large disturbances in water waves are described by Chamberlain 
and Porter \cite{kn:Chamberlain}. 

It is well known that the spectra of non-selfadjoint operators can be 
pathologically sensitive to perturbation of the operator. Matrix examples
of such sensitivity are provided in the classic text of Wilkinson \cite{kn:Wilkinson}. 
For a recent study in the context of non-selfadjoint Sturm-Liouville operators 
see Davies \cite{kn:AAD}, and for a study in the context
of general operators via {\em pseudospectra} see Trefethen \cite{kn:Trefethen}.  
Given this sensitivity, it seems important to ask: under what conditions can
one expect the regularization process used for selfadjoint singular
Sturm-Liouville operators to be successful for non-selfadjoint Sturm-Liouville
operators? In particular, can one recover results on spectral inclusion and
spectral exactness? If not, might one at least be able to recover results
on pseudospectral inclusion and pseudospectral exactness, or develop
a-posteriori tests for spectral exactness?

We seek to answer these questions in this paper, for singular second order
non-selfadjoint Sturm-Liouville problems. 

For selfadjoint problems a singular endpoint is either of {\em limit point }
or of {\em limit circle} type. This is the Titchmarsh-Weyl theory and may
be developed either using methods of complex analysis (see Titchmarsh 
\cite{kn:Titchmarsh}) or using the theory of deficiency indices for symmetric
operators on Hilbert spaces (see, e.g., Dunford and Schwartz \cite{kn:Dunford}).
The analogous theory for non-selfadjoint problems is due to Sims \cite{kn:Sims}
and to Brown, Evans, McCormack and Plum \cite{kn:BEMP}. It is based on the
Titchmarsh approach to the selfadjoint case, and will be very important in
this paper. The other ingredient which we shall find useful is a
non-selfadjoint analogue of the results of Reed and Simon on spectral inclusion
and spectral exactness \cite[theorems VIII.23-VIII.25]{kn:RS}, for which we shall 
use results from Harrabi \cite{kn:Harrabi} and Kato \cite[p. 208]{kn:Kato}.

\section{A review of the Sims Classification}\label{section:2}
The problem which we consider concerns the spectral behavior of
 \beq
{\cal M}[y]= \frac{1}{w} [ - (py^{'})^{'} + q y ] \;\;\;{\rm on}\;\;[a,b),
\label{eq:2.1m}
\enq
where as usual
\newcounter{rem1}
\begin{list}%
{(\roman{rem1})}{\usecounter{rem1}
\setlength{\rightmargin}{\leftmargin}}
 \item
$w>0$, $p \neq 0$ a.e. on $[a,b)$ and $w, 1/p \in L^1_{loc}[a,b)$;
\item
$p,q$ are complex-valued, $q \in L^1_{loc}[a,b)$ and
\beq
Q=\overline{co} \left\{ \frac{q(x)}{w(x)} + r  p(x):
 x \in [a,b), \; 0 < r < \infty \right\} \neq \Cc. \label{eq:2.2m}
\enq
\end{list}
These assumptions imply that $a$ is a regular point of (\ref{eq:2.1m})
and  we shall assume that $b$ is a singular point. By this we mean that either $b=
+\infty$ or that $\int_a^b(w +\frac{1}{\mid p \mid} + \mid q \mid ) dx = \infty$.
Since we are assuming that $Q$ does not occupy all of $\Cc$, it is known that its complement 
has either one or two connected components.
For $\lambda_0 \in \Cc \backslash Q $ we denote by $K=K(\lambda_0)$ the nearest point
in $Q$ to $\lambda_0$ and by $L$ the tangent to $Q$ at $K$ and  arrange by translation and rotation 
through an angle $\eta$ for $L$  to coincide with the imaginary axis
while $\lambda_0$ and $Q$ are contained in the new left and right half planes respectively.
That is, for all $x \in [a,b)$ and $r \in (0,\infty)$, we require by choice of $K$ and $\eta$
that
\beq
\Re[ \{ r p(x) + \frac{ q(x)}{w(x)}-K \} e^{i \eta} ] \geq 0 \label{eq:2.3m}
\enq
and
\beq
\Re [ (\lambda_0-K) e^{i \eta} ] <0. \label{eq:2.3am}
\enq
   The set of all such {\it admissible } pairs $(\eta,K)$ we call $S$ and we also define
\[ \lek= \{ \lambda \in {\Cc} : \Re [(\lambda-K) e^{i \eta}]\leq 0 \}. \]
In order to obtain from (\ref{eq:2.1m}) a well posed eigenvalue problem
we need to introduce boundary conditions at $a$ and possibly at $b$.  
The conditions at $a$ will be given in the form
\be y(a)\cos\alpha + py'(a)\sin\alpha = 0, \label{eq:bca} \ee
where the parameter $\alpha$, which may be complex, will be subject to the  condition 
\beq
\Re[ e^{i \eta} \cos \alpha \;\overline{ \sin \alpha}] \leq 0.
\label{eq:2.6m} 
\enq
This gives rise to a set $S(\alpha)$  which is defined as the subset of $S$ in which  
(\ref{eq:2.6m}) holds.   We note that $\alpha=0$ and $\alpha=\pi/2$ correspond to Dirichlet 
and Neumann boundary conditions respectively. 
\par
When  $p=1$ and  $q$ is real the classical theory of Weyl \cite{Weyl} and Titchmarsh 
\cite{kn:Titchmarsh} shows that if $\theta$ and $\phi$ are linearly independent solutions of 
(\ref{eq:2.1m}) which satisfy
\begin{eqnarray}
\phi(a,\lambda)=\sin \alpha, & \theta(a,\lambda)= \cos \alpha, \nonumber  \\
p\phi^{'}(a,\lambda)=-\cos \alpha, & p\theta^{'}(a,\lambda)= \sin \alpha,\label{eq:2.6a}
\end{eqnarray}
where $\alpha$ is now real, then there is a complex number  $m(\lambda)$, a function of the 
strictly complex variable $\lambda$, such that 
\be
\psi=\theta+m\phi
\label{eq:mdef} \ee
lies in $L^2_w[a,b)$. When, up to constant multiples, $\psi$ is the only solution of the differential
equation which lies in $L^2_w[a,b)$, we say that (\ref{eq:2.1m})
is in the limit point case at $b$.  If however both $\theta$ and  $\phi$
lie in $L^2_w[a,b)$ then we say that (\ref{eq:2.1m}) is in the limit circle case at $b$. 
In this case an additional boundary condition at $b$ is needed in order to
make (\ref{eq:2.1m},\ref{eq:bca}) into a well-posed eigenvalue problem. There is
a one to one correspondence between this additional boundary condition and the
choice of function $m(\cdot)$ in (\ref{eq:mdef}), in the sense that with an
appropriate choice of boundary condition there exists a unique function
$m(\cdot)$ such that eqn. (\ref{eq:mdef}) defines a solution $\psi$
of the differential equation satisfying the boundary condition at $x=b$, while
with an allowed choice of $m(\cdot)$ the function $\psi$ defined in
(\ref{eq:mdef}) can itself be used, for appropriate $\lambda$, to define the
boundary condition at $x=b$ in the form $[y,\psi](b)=0$, where $[f,g]:=p(fg'-f'g)$ 
denotes the Wronsian of two functions $f$ and $g$. 
It is known that the classification of limit point or limit circle is independent of the 
strictly complex parameter $\lambda$. The terminology of limit point or limit circle owes 
its origin  to the method used to establish the existence and possible uniqueness of $\psi$
in (\ref{eq:mdef}). It may be shown that the spectral points of any realisation of 
(\ref{eq:2.1m}) as an operator in $L^2_w[a,b)$ may be charecterised by the behaviour  
in the limit as $\Im\lambda \rightarrow 0$ of the function $m(\lambda)$ associated with the 
boundary conditions defining the domain of the realisation.
\par
Many of these notions may be carried over to the case when $p$, $q$ and $\alpha$ are complex.
In a seminal paper Sims \cite{kn:Sims} shows that when $p=w=1$ and $\Im q \in \Cc_{-}$, 
where $\Cc_{-}$ denotes the strictly lower complex plane,
then the limit point / limit circle classification of Weyl now gets replaced by a threefold classification.
We shall discuss this in the more general setting of \cite{kn:BEMP} which only requires 
(\ref{eq:2.2m}), (\ref{eq:2.3m}) and (\ref{eq:2.3am}) to hold. Using a nesting circle  
method based on that of both Weyl and Sims, Brown et al. prove  the following theorem.
\begin{theorem} \label{theorem:b1} \cite{kn:BEMP}
For $\lambda \in \lek$, $(\eta,K) \in S(\alpha)$ the following distinct cases are possible, 
the first two being sub-cases of the limit point case:
\begin{itemize}
\item
Case I : there exists a unique solution of (\ref{eq:2.1m}) satisfying 
\beq
\int_a^b\Re[ e^{i \eta}  \{ p \mid y' \mid^2 +(q-K w ) \mid y \mid^2 \}]dx 
 + \int_a^b \mid y \mid^2 w dx < \infty;
\label{eq:2.25m}
\enq
 and this is the only solution satisfying $y \in L^2_w[a,b)$;
\item
Case II : there exists a unique solution of (\ref{eq:2.1m}) satisfying (\ref{eq:2.25m}), but all 
solutions of (\ref{eq:2.1m}) lie in $L^2_w[a,b)$;
\item
Case III: all solutions of (\ref{eq:2.1m}) lie in $L^2_w[a,b)$ and satisfy
(\ref{eq:2.25m}).
\end{itemize}
It may also be shown that the classification is independent of $\lambda$ in the sense that
\begin{list}%
{(\roman{rem1})}{\usecounter{rem1}
\setlength{\rightmargin}{\leftmargin}}
 \item
if all solutions of (\ref{eq:2.1m}) satisfy  (\ref{eq:2.25m}) for some $\lambda^{'} \in \lek$ (i.e. Case III) 
then all solutions of (\ref{eq:2.1m}) satisfy (\ref{eq:2.25m}) for all $\lambda \in \Cc$;
\item
if  all solutions of (\ref{eq:2.1m})  lie in $L^2_w[a,b)$ for some $ \lambda^{'} \in \Cc$ then 
all solutions of (\ref{eq:2.1m}) satisfy $y \in L^2_w[a,b)$ for all $ \lambda \in \Cc$.
\end{list}
\end{theorem}
It is interesting to examine the case when $p$ is real and non-negative.
In this case   for some 
$\eta \in [ -\frac{\pi}{2},\frac{\pi}{2}]$ and $K \in \Cc$  let
\beq
\theta_{K,\eta}(x) = \Re [ e^{i \eta}( q(x)-Kw(x))] \geq 0\;\; a.e. \; x \in(a,b). 
\label{eq:2.27}
\enq
Then the condition (\ref{eq:2.25m}) in the Sims characterisation of (\ref{eq:2.1m}) 
in Theorem \ref{theorem:b1} for   $ \lambda \in \lek$, $(\eta,K) \in S(\alpha)$, becomes
\beq
\cos \eta \int_a^b p\mid y' \mid^2dx + \int_a^b \theta_{K\eta}(x) \mid y(x) \mid^2 dx +   
\int_\alpha ^b \mid y(x)\mid^2w(x) dx < \infty.
\label{eq:2.28m}
\enq
 In this case the remark on the independence of the classification  can  be extended to the following:
\begin{list}%
{(\roman{rem1})}{\usecounter{rem1}
\setlength{\rightmargin}{\leftmargin}}
 \item
if for some $\lambda^{'} \in \Cc$ all the solutions of (\ref{eq:2.1m}) satisfy
(\ref{eq:2.28m}), then for all $\lambda \in \Cc$ all solutions of (\ref{eq:2.1m}) 
satisfy (\ref{eq:2.28m});
\item
if  for some $\lambda^{'} \in \Cc$ all the solutions of (\ref{eq:2.1m}) satisfy one of
 \beq
\cos \eta \int_a^b p \mid y' \mid^2 dx < \infty, \label{eq:2.29m}
\enq
\beq
\int_a^b \theta_{K\eta }\mid y\mid^2dx < \infty, \label{eq:2.30m}
\enq
then the same applies for all $\lambda \in \Cc$.
\end{list}
We remark that Sim's analysis is the special case of the above when $\eta=\pi/2$, $K=0$. This restriction 
overlooks the interesting features present in (\ref{eq:2.28m}) when $\eta \in ( -\frac{\pi}{2},\frac{\pi}{2})$, 
namely, that the classification in Theorem \ref{theorem:b1} involves a weighted Sobolev space as well as 
$L^2_w[a,b)$. The paper \cite{kn:BEMP} also examines the analytic behaviour of $m(\lambda)$ and the connection 
between  this and the spectrum of $M$, an operator realisation of ${\cal M}$ in $L^2_w[a,b)$. This 
is sumarised in the following theorems, in which $m$ denotes the unique function such that (\ref{eq:mdef}) 
defines a solution of (\ref{eq:2.1m}) which either (i) lies in $L^2_w[a,b)$ (Sims Case I) or (ii)
lies in $L^2_w[a,b)$ and satisfies the additional boundary condition at $x=b$ (Sims
Cases II and III).

\begin{theorem} (Theorem 4.7, \cite{kn:BEMP}) \label{theorem:b2}
In Cases II and III, $\lambda_0$ is a pole of $m$ of order $s$ if and only if $\lambda_0$ is an eigenvalue of 
$M$ of algebraic multiplicity $s$.
\end{theorem} 
\begin{theorem} (Theorem 4.13, \cite{kn:BEMP}) \label{theorem:b3}
Suppose that (\ref{eq:2.2m}) is in Case I. Define
\[ Q(\alpha) = \bigcap_{(\eta,K)\in S(\alpha)} (\Cc \backslash \Lambda_{\eta,K}), \]
and let $Q_c(\alpha)$ denote the set $Q(\alpha)$ when the underlying interval 
is $[c,b)$ rather than $[a,b)$. Define
\begin{eqnarray*}
Q_c &:= &\overline{co} \left\{ \frac{ q(x)}{w(x)}+rp(x) : x \in [c,b),  \; r \in (0,\infty)
\right\}, \label{eq:3.6} \\
Q_b &:=& \cap_{ c \in (a,b)} Q_c , \;\;\; Q_b(\alpha) = \cap_{c \in (a,b) } Q_c(\alpha),  
\end{eqnarray*}
Then $m(\lambda)$ is defined throughout $\Cc \backslash Q (\alpha)$ and has a meromorphic 
extension to $\Cc \backslash Q_b(\alpha)$, with poles only in 
$  Q(\alpha) \backslash Q_b (\alpha)$. In addition  $\lambda$ is a pole of $m(\lambda)$ if and only if
 $\lambda$ is an eigenvalue of $M$ for $\lambda \not \in Q_b(\alpha)$.
\end{theorem}

\section{Tests for spectral inclusion and spectral exactness}
In this section we prove a simple theorem (Theorem \ref{theorem:4.1}) which 
allows us to test a convergent sequence of eigenvalue approximations obtained from 
a sequence of truncated  interval problems, in order to determine whether or not 
the limit of the sequence is truly an eigenvalue of our original problem. We also 
prove two additional results (Theorem \ref{theorem:4.2} and Theorem \ref{theorem:4.3})
which give methods for determining whether or not the hypotheses
of Theorem \ref{theorem:4.1} are satisfied for a given problem. The second of these, 
Theorem \ref{theorem:4.3}, extends a convergence result in 
\cite{kn:BEMP} from the complement of the numerical range of our singular 
operator into a set which is typically much larger. 

Finally, we show how Theorem \ref{theorem:4.2} and Theorem \ref{theorem:4.3}
allow us to develop a test for spectral inclusion, to ensure that there will
be no eigenvalues of the original problem which remain unapproximated by
the truncation process.

\subsection{Spectral Exactness}

We denote by $m(\cdot)$ the $m$-function developed in section \ref{section:2};
in Sims cases II and III, this corresponds to a particular choice of boundary
condition at $x=b$. From Theorems \ref{theorem:b2} and \ref{theorem:b3} 
we know that the poles of $m$ are the eigenvalues of a realization $M$
of the differential operator ${\cal M}$ subject to a boundary condition of the form
\be (\cos \alpha) y(a) + (\sin \alpha) py'(a) = 0 \label{eq:4.1} \ee
(and possibly an additional boundary condition at $x=b$).

We denote by $L$ the operator, and by $\ell(\cdot)$ the Titchmarsh-Weyl function,
when the boundary condition at $x=a$ is changed to
\be (\sin \alpha) y(a) - (\cos \alpha) py'(a) = 0, \label{eq:4.2} \ee
but the boundary conditions at $x=b$ (where applicable) are left unchanged,
the same as those for $M$.

The functions $\ell$ and $m$ are related by the identity
\[ m(\lambda)\ell(\lambda) = -1 \]
(see \cite[eqn. 5.17]{kn:BEMP}).

Let $(b_n)_{n\in \Nn}$ be a sequence such that $b_n\nearrow b$ as $n\nearrow 
\infty$. Following \cite[\S 2]{kn:BEMP} we may construct a sequence $M_n$ of regular 
operators defined on the intervals $[a,b_n]$. These operators $M_n$ are
still defined by $M_n y = {\cal M} y$ on their domains: it is the boundary
conditions defining these domains which are of interest. At $x=a$ we keep
the boundary condition (\ref{eq:4.1}). At $x=b_n$ we impose a boundary condition
of the form
\be y(b_n)\cos\beta_n + py'(b_n)\sin\beta_n = 0. \label{eq:bcnew} \ee
The Titchmarsh-Weyl function $m_n$ associated with 
$M_n$ is then given in terms of the solutions $\theta$ and $\phi$ of 
(\ref{eq:2.6a}) by
\[ m_n(\lambda) = -\frac{\theta(b_n,\lambda)\cot\beta_n 
 + p\theta'(b_n,\lambda)}{\phi(b_n,\lambda)\cot\beta_n + p\phi'(b_n,\lambda)}. \]
\cite[eq. (2.6)]{kn:BEMP}. We shall examine in Lemma \ref{lemma:extra},
Theorem \ref{theorem:4.2} and Theorem \ref{theorem:4.3} below conditions
on the $\beta_n$ which ensure that, for $\lambda$ in certain regions of
$\Cc$,
\be m(\lambda) = \lim_{n\rightarrow\infty}m_n(\lambda). \label{eq:4.3} \ee
Let $L_n$ be a sequence of regular operators defined on the
intervals $[a,b_n]$ with boundary condition (\ref{eq:4.2}) at $x=a$ and with the
same boundary conditions (\ref{eq:bcnew}) as the $M_n$ at $x=b_n$, so that the
associated Titchmarsh-Weyl functions $\ell_n(\cdot)$ satisfy
\[ m_n(\lambda)\ell_n(\lambda) = -1. \]
By analogy with (\ref{eq:4.3}) we shall assume that in some appropriate regions
of $\Cc$, 
\be \ell(\lambda) = \lim_{n\rightarrow\infty}\ell_n(\lambda). \label{eq:4.4} \ee

\begin{theorem}\label{theorem:4.1} (Test for spectral inexactness) \,
Suppose that $\mu\in \Cc$ has the following properties:
\begin{enumerate}
\item $\mu$ does not lie in the spectrum of $M$;
\item $m_n\rightarrow m$ and $\ell_n\rightarrow \ell$ uniformly on any compact
 annulus of sufficiently small outer radius surrounding $\mu$.
\end{enumerate}
Then there are only two possibilities:
\begin{description}
\item[(a)] there exists a neighbourhood ${\cal N}$ of $\mu$ and $N\in \Nn$ such that no
 eigenvalue of $M_n$ lies in ${\cal N}$ for any $n\geq N$;
\item[(b)] there exists a monotone increasing sequence $(n_j)_{j\in \Nn}$ of
 positive integers, and two associated sequences $(\lambda_j)_{j\in \Nn}$,
 $(\mu_j)_{j\in \Nn}$, such that 
\[ \lambda_j \in \sigma(M_{n_j}), \;\;\;
   \mu_j \in \sigma(L_{n_j}), \]
and $\lambda_j \rightarrow \mu$, $\mu_j \rightarrow \mu$ as $j\rightarrow
\infty$.
\end{description}
\end{theorem}
A consequence of this theorem is that if a subsequence of eigenvalues of 
the regularized operators $M_n$ converges to some point $\mu$ which is not
an eigenvalue (spectral inexactness) then the $L_n$ will also possess
a subsequence of eigenvalues converging to the same point $\mu$. Moreover, if 
subsequences of eigenvalues of $M_n$ and of $L_n$ converge to the same point 
then at least one of the subsequences is spectrally inexact, because the 
boundary conditions (\ref{eq:4.1}) and (\ref{eq:4.2}) ensure that $M$ and 
$L$ have no shared eigenvalues. Theorem \ref{theorem:4.1} therefore gives 
us a test for spectral exactness: if only the $M_n$, and not the
$L_n$, possess eigenvalues accumulating at $\mu$,
then $\mu$ must be an eigenvalue of $M$.
\vspace{2mm}

\noindent {\bf Proof} of Theorem \ref{theorem:4.1}. Let $A$ be any 
sufficiently small annulus surrounding $\mu$ and let $\Gamma$ be a closed
contour in $A$ surrounding $\mu$. For any function $f$ which is meromorphic in 
a simply connected open set containing $\Gamma$ we denote by $N_Z(f,\Gamma)$ the
number of zeros of $f$ inside $\Gamma$ and by $N_P(f,\Gamma)$ the
number of poles of $f$ inside $\Gamma$. Rouch\'{e}'s Theorem gives
\[ N_Z(m,\Gamma) - N_P(m,\Gamma) = \frac{1}{2\pi i}\int_\Gamma
\frac{m'(\lambda)}{m(\lambda)}d\lambda. \]
In view of the identity $m(\lambda)\ell(\lambda) = -1$ we have
$N_Z(m,\Gamma) = N_P(\ell,\Gamma)$ and $1/m = -\ell$, so we
can write this
\be N_P(\ell,\Gamma) - N_P(m,\Gamma) = -\frac{1}{2\pi i}\int_\Gamma
 m'(\lambda)\ell(\lambda)d\lambda. \label{eq:4.5}\ee
As $\mu$ does not lie in the spectrum of $M$, we have $N_P(m,\Gamma)=0$
for all sufficiently small annuli $A$. Hence 
\be N_P(\ell,\Gamma)
 = -\frac{1}{2\pi i}\int_\Gamma m'(\lambda)\ell(\lambda)d\lambda. \label{eq:4.6} \ee
By arguments similar to those which gave (\ref{eq:4.5}) we have
\[ N_P(\ell_n,\Gamma) - N_P(m_n,\Gamma) = 
-\frac{1}{2\pi i}\int_\Gamma m_n'(\lambda)\ell_n(\lambda)d\lambda. \]
The uniform convergence $m_n\rightarrow m$ implies uniform convergence
of $m_n'$ to $m'$ (by the Cauchy integral representation of the derivative).
Combined with the uniform convergence $\ell_n\rightarrow \ell$
this yields
\[ N_P(\ell_n,\Gamma) - N_P(m_n,\Gamma) = 
-\frac{1}{2\pi i}\int_\Gamma m'(\lambda)\ell(\lambda)d\lambda \]
for all sufficiently large $n$. Combining this with (\ref{eq:4.6}) 
we have, for all sufficiently large $n$,
\be N_P(\ell,\Gamma) = N_P(\ell_n,\Gamma) - N_P(m_n,\Gamma). 
\label{eq:4.7} 
\ee
In the case of possibility {\bf (a)}, the sequence $M_n$ does not have any 
eigenvalues converging spuriously to $\mu$, which is a non-eigenvalue of $M$.
Since eigenvalues of $M_n$ are poles of $m_n$ we have $N_P(m_n,\Gamma) = 0$
and (\ref{eq:4.7}) then shows that the sequence $L_n$ is spectrally
exact for $L$ near $\mu$. Thus we have spectral exactness near $\mu$ for both 
$M$ and $L$. In the case that {\bf (a)} is not true, then for some arbitrarily small
annuli $A$ and arbitrarily large $n\in \Nn$ we will have
\[ N_P(m_n,\Gamma) = \mbox{no. of eigenvalues of $M_n$ inside $\Gamma$} > 0, \]
and, from (\ref{eq:4.7}),
\[ N_P(\ell_n,\Gamma) - N_P(m_n,\Gamma) = N_P(\ell,\Gamma) \geq 0, \]
which shows that $N_P(\ell_n,\Gamma)$, the number of eigenvalues of $L_n$
inside $\Gamma$, is at least 1. Thus $L_n$ also has an eigenvalue close
to $\mu$. This gives possibility {\bf (b)}. \hfill $\Box$

\subsection{Conditions for $m$-function convergence}

We now examine the hypotheses $m_n\rightarrow m$ and $\ell_n\rightarrow\ell$ 
of Theorem \ref{theorem:4.1}. Under what conditions do these hold? 

We consider first Sims Cases II and III.  The following result explains
how to choose the boundary condition (\ref{eq:bcnew}) to ensure that
(\ref{eq:4.3}) and (\ref{eq:4.4}) hold.
\begin{lemma}\label{lemma:extra} 
Suppose that the differential equation is of Sims Case II or III.
Let $\lambda' \in \Lambda_{\eta,K}$ be fixed.
Express the boundary conditions at $x=b$ for $M$ in terms of an $L^2_w$-solution 
$\psi(\cdot,\lambda') = \theta(\cdot,\lambda') + m(\lambda')\phi(\cdot,\lambda')$ 
of the differential equation ${\cal M}\psi = \lambda' \psi$, in the form
\[ [y,\psi(\cdot,\lambda')](b) = 0, \]
\cite[eq. (4.11)]{kn:BEMP} where $[\cdot,\cdot]$ denotes the usual 
Wronskian $[u,v] := p(uv'-u'v)$. Then appropriate boundary conditions 
(\ref{eq:bcnew}) are given by choosing
\be (\cos\beta_n,\sin\beta_n) =
const.(p\psi'(b_n,\lambda'),-\psi(b_n,\lambda')) \label{eq:tbc} \ee
so that (\ref{eq:bcnew}) is simply the condition $[y,\psi](b_n) = 0$.
\end{lemma}
\noindent {\bf Proof} \, Define  $\psi_n(\cdot,\lambda') = \theta(\cdot,\lambda') + 
m_n(\lambda')\phi(\cdot,\lambda')$, in which $m_n$ is chosen so that $\psi_n$
satisfies (\ref{eq:bcnew}) with $\beta_n$ given by (\ref{eq:tbc}). From the
definition of the $\beta_n$, the fact that $\psi_n$ satisfies (\ref{eq:bcnew})
may be written as 
\[ [\psi_n(\cdot,\lambda'),\psi(\cdot,\lambda')](b_n) = 0. \]
Substituting $\psi_n(\cdot,\lambda') = \theta(\cdot,\lambda') +
m_n(\lambda')\phi(\cdot,\lambda')$ into this equation yields
\be m_n(\lambda') = -\frac{[\theta(\cdot,\lambda'),\psi(\cdot,\lambda')](b_n)}{[\phi
(\cdot,\lambda'),\psi(\cdot,\lambda')](b_n)}. \label{eq:mnew1} \ee
In the identity $[\psi(\cdot,\lambda'),\psi(\cdot,\lambda')](b) = 0$,
replace the first instance of $\psi(\cdot,\lambda')$ by
$\theta(\cdot,\lambda') + m(\lambda')\phi(\cdot,\lambda')$, and hence obtain
\be m(\lambda') = -\frac{[\theta(\cdot,\lambda'),\psi(\cdot,\lambda')](b)}{[\phi
(\cdot,\lambda'),\psi(\cdot,\lambda')](b)}. \label{eq:mnew2} \ee
Comparing (\ref{eq:mnew1}) with (\ref{eq:mnew2}) establishes (\ref{eq:4.3})
when $\lambda = \lambda'$:
\be \lim_{n\rightarrow\infty}m_n(\lambda') = m(\lambda'). \label{eq:mnew3} \ee
Now Brown et al. \cite[Corollary 3.4, eqn. (3.4)]{kn:BEMP} give a formula which 
allows us to extend this result to other values of $\lambda$:
\be m(\lambda) =
 \frac{m(\lambda') - (\lambda-\lambda')\int_a^b w(x) 
\theta(x,\lambda)\psi(x,\lambda')dx}{1+(\lambda-\lambda')\int_a^b w(x)
\phi(x,\lambda)\psi(x,\lambda')dx  }. \label{eq:4.8} 
\ee
This formula possesses the regular-interval analogue
\be m_n(\lambda) =
 \frac{m_n(\lambda') - (\lambda-\lambda')\int_a^{b_n} w(x) 
\theta(x,\lambda)\psi_n(x,\lambda')dx}{1+(\lambda-\lambda')\int_a^{b_n} w(x)
\phi(x,\lambda)\psi_n(x,\lambda')dx  }. \label{eq:4.9} 
\ee
These formulae hold at any point which is not an eigenvalue of $M$ or of
$M_n$, respectively. Moreover, since the equation is in Sims Case II or III,  
all of $\theta(\cdot,\lambda)$, $\theta(\cdot,\lambda')$, $\phi(\cdot,\lambda)$
and $\phi(\cdot,\lambda')$ lie in $L^2_w[a,b)$. Using
$\psi_n(\cdot,\lambda') - \psi(\cdot,\lambda') = (m_n(\lambda')-
m(\lambda'))\phi(\cdot,\lambda')$ it follows from (\ref{eq:mnew3})
that $\psi_n(\cdot,\lambda')\rightarrow \psi(\cdot,\lambda')$ in
$L^2_w[a,b)$. Hence, combining (\ref{eq:4.8}) and
(\ref{eq:4.9}), we obtain the convergence
\[ \lim_{n\rightarrow\infty}m_n(\lambda) = m(\lambda) \]
at any point $\lambda$ which is not an eigenvalue of $M$. This establishes
(\ref{eq:4.3}), and (\ref{eq:4.4}) is proved similarly. \hfill $\Box$
\vspace{2mm}

Of course, Theorem \ref{theorem:4.1} requires more than just pointwise
convergence, and so it is fortunate that the following
stronger result holds.
\begin{theorem}\label{theorem:4.2} Suppose that the problem is Sims Case II
or Sims Case III at $x=b$ and let the hypotheses of Lemma \ref{lemma:extra}
hold. Then for $\lambda$ in any compact set 
${\cal K}\subseteq \Cc$ not containing eigenvalues of $M$, 
\be \lim_{n\rightarrow\infty}m_n(\lambda) = m(\lambda), \label{eq:4.10} \ee
the convergence being uniform over ${\cal K}$. 
\end{theorem}
\noindent {\bf Proof} \, That $m_n(\lambda)\rightarrow m(\lambda)$ pointwise
on ${\cal K}$ has already been proved in Lemma \ref{lemma:extra}. The uniformity of the 
convergence depends on having a uniform bound on the $L^2_w$ norms of $\theta(\cdot,\lambda)$
and $\phi(\cdot,\lambda)$ for $\lambda\in {\cal K}$. This can be obtained by a standard 
variation of parameters argument, expressing the solutions in terms of 
$\theta(\cdot,\lambda')$ and $\phi(\cdot,\lambda')$: see Sims 
\cite[section 3, Theorem 2]{kn:Sims} and also \cite[Remark 2.2]{kn:BEMP}.
\hfill $\Box$
\vspace{2mm}

For the functions $\ell_n$ and $\ell$, a result exactly analogous to Theorem
\ref{theorem:4.2} is clearly valid, the only difference being now that ${\cal K}
$ must not contain eigenvalues of $L$. 
\vspace{2mm}

We turn now to Sims Case I. In order to handle this case it
is necessary to know more about the behaviour of the solutions of the
differential equation. Suppose that for some $\lambda\in\Cc$, the differential
equation possesses `small' and `large' solutions. We shall assume that 
the small solution is the (unique up to scalar multiples) square integrable
solution $\psi(x,\lambda)$, and we denote the non-unique large solution 
by $\Upsilon(x,\lambda)$. By `small' and `large' we mean that 
these solutions satisfy the condition
\be \lim_{x\rightarrow b} \frac{\psi(x,\lambda)}{\Upsilon(x,\lambda)}
 = 0. \label{eq:4.11} 
\ee
Clearly $\Upsilon$ is not unique: $\Upsilon + \psi$, for example, is also
a `large' solution in the sense of (\ref{eq:4.11}).  The solutions $\theta$
and $\phi$ of (\ref{eq:2.6a}) can clearly be written in terms of
$\psi$ and $\Upsilon$:
\be \begin{array}{c}
 \theta(x,\lambda) = c_1\psi(x,\lambda) + c_2\Upsilon(x,\lambda), \\
 \phi(x,\lambda) = d_1\psi(x,\lambda) + d_2\Upsilon(x,\lambda), 
\end{array} \label{eq:4.14}
\ee
in which the constants $c_1$, $c_2$, $d_1$ and $d_2$ are given by
\be c_1 =  \left((\cos\alpha) p\Upsilon'(a,\lambda)
                 -(\sin\alpha) \Upsilon(a,\lambda)\right)/W, \;\;\;
    c_2 =  \left(-(\cos\alpha) p\psi'(a,\lambda)
                 +(\sin\alpha) \psi(a,\lambda)\right)/W, 
\label{eq:4.16}\ee 
\be d_1 =  \left((\sin\alpha) p\Upsilon'(a,\lambda)
                 +(\cos\alpha) \Upsilon(a,\lambda)\right)/W, \;\;\;
    d_2 =  \left(-(\sin\alpha) p\psi'(a,\lambda)
                 -(\cos\alpha) \psi(a,\lambda)\right)/W, 
\label{eq:4.17}\ee 
where $W = p(\psi\Upsilon'-\psi'\Upsilon)$ is the usual Wronskian.
Suppose that $m_n$ is defined by the requirment that the solution
\[ \psi_n(\cdot,\lambda) = \theta(\cdot,\lambda) + m_n(\lambda)
\phi(\cdot,\lambda) \]
satisfy the boundary condition $\psi_n(b_n,\lambda) = 0$. Then 
\be m_n(\lambda) = -\frac{\theta(b_n,\lambda)}{\phi(b_n,\lambda)}. 
\label{eq:4.15} \ee
Now combining (\ref{eq:4.11}) with (\ref{eq:4.14}) we have
\[ \theta(b_n,\lambda) \sim c_2\Upsilon(b_n,\lambda), \;\;\;
   \phi(b_n,\lambda) \sim d_2 \Upsilon(b_n,\lambda) \]
for large $n$. Combining this with (\ref{eq:4.15}) yields
\[ m_n(\lambda) \sim -\frac{c_2}{d_2} \]
for large $n$. Together with (\ref{eq:4.16}) and (\ref{eq:4.17})
this yields, for large $n$,
\be m_n(\lambda) \sim 
\frac{-(\cos\alpha) p\psi'(a,\lambda)+(\sin\alpha)\psi(a,\lambda)}{
(\sin\alpha) p\psi'(a,\lambda)+(\cos\alpha) \psi(a,\lambda)}
 = m(\lambda), \label{eq:4.18} \ee
the last equality in (\ref{eq:4.18}) being an immediate consequence
of \cite[Definition 4.10]{kn:BEMP}. From these considerations the following
result is clearly true.

\begin{theorem} \label{theorem:4.3}
Suppose that the differential equation is of Sims Case I type at $x=b$.
Let ${\cal K}\subseteq \Cc$ be a compact set such that for $\lambda\in 
{\cal K}$ the square-integrable solution $\psi(x,\lambda)$ of the differential 
equation exists and is an analytic function of $\lambda$. Suppose moreover that
there exists a second solution $\Upsilon(x,\lambda)$ such that
\be \lim_{x\rightarrow b} \frac{\psi(x,\lambda)}{\Upsilon(x,\lambda)}
 = 0, \label{eq:4.19} \ee
the limit being uniform with respect to $\lambda \in {\cal K}$.  Suppose that
the domains of the operators $M_n$ are determined by the Dirichlet conditions 
$y(b_n)=0$. Then we have the convergence
\[ \lim_{n\rightarrow\infty}m_n(\lambda) = m(\lambda) \]
uniformly for $\lambda\in {\cal K}$.
\end{theorem}
Clearly a similar result holds for the functions $\ell_n(\lambda)$ and their
convergence to the function $\ell(\lambda)$.
\vspace{2mm}

\noindent {\bf Remark} The result of Theorem \ref{theorem:4.3} will also
hold if the domains of the $M_n$ are defined by certain other boundary conditions
at $x=b_n$. Suppose that a boundary condition
\[ y(b_n)\cos\beta_n  + y'(b_n)\sin\beta_n = 0 \]
is imposed, where the $\beta_n$ are complex numbers. Then it may be shown
that the result continues to hold provided 
\be \lim_{n\rightarrow\infty}\frac{\psi(b_n,\lambda)\cos\beta_n + 
         \psi'(b_n,\lambda)\sin\beta_n}{\Upsilon(b_n,\lambda)\cos\beta_n 
         + \Upsilon'(b_n,\lambda)\sin\beta_n} = 0, \label{eq:abc} \ee
locally uniformly with respect to $\lambda$.
In problems where $\psi'(x,\lambda)/\Upsilon'(x,\lambda)\rightarrow 0$
as $x\rightarrow b$, one would have to choose the $\beta_n$ quite
carefully for (\ref{eq:abc}) to fail.

\subsection{A simple test for spectral inclusion}
In the selfadjoint case, spectral inclusion is usually very easy to prove. In
fact, suppose $T$ is a selfadjoint operator on a domain $D(T)$ in a Hilbert 
space $H$ and let $(T_n)$ be a sequence of operators with domains $(D(T_n))$ which 
converge pointwise to $T$ on some set 
${\cal C}\subseteq \cap_{m\in\Nn}\cup_{n\geq m}D(T_n)$:
\be \lim_{n\rightarrow\infty}\| T_n f - Tf \| = 0 \;\;\;
 \forall f\in {\cal C}. \label{eq:4.20} \ee
Then provided ${\cal C}$ is a core of $T$ -- in other words, provided
the set $({\cal C}, T{\cal C})$ is dense in the graph of $T$ -- the sequence
$(T_n)$ will be spectrally inclusive for $T$: every eigenvalue of $T$ will
be the limit of some sequence $(\lambda^{(n)})$ in which $\lambda^{(n)}$
lies in the spectrum of $T_n$.

In the non-selfadjoint case a result of such generality does not seem to exist,
although some of the results of Harrabi \cite{kn:Harrabi} come quite close.
We shall examine some corollaries of Harrabi's work, as well as a standard
result from Kato \cite{kn:Kato}, in section \ref{subsection:4.4} below. In 
this subsection, however, we shall show that spectral inclusion always holds 
in Sims Cases II and III, and in Sims Case I in those parts of the complex 
plane where Theorem \ref{theorem:4.3} holds.

\begin{theorem}\label{theorem:4.4} (Test for Spectral Inclusion)\,
In the notation of Theorem \ref{theorem:4.1}, suppose that $\mu\in \Cc$ is an 
isolated eigenvalue of $M$. Suppose also that $m_n(\lambda)\rightarrow m(
\lambda)$ as $n\rightarrow\infty$ uniformly on any compact annulus 
surrounding $\mu$. Then there exists a sequence $(\lambda^{(n)})$ in which
$\lambda^{(n)}$ lies in the spectrum of $M_n$, such that
$\lim_{n\rightarrow\infty}\lambda^{(n)} = \mu$.
\end{theorem}
\noindent {\bf Proof}\, Given $\epsilon>0$ sufficiently small, surround
$\mu$ by an annulus $A$ whose outer radius is at most $\epsilon$ and let
$\Gamma$ be a circular contour surrounding $\mu$ and contained in $A$.
Since $\mu$ is a pole of $m$, we may assume by taking $\epsilon$ sufficiently
small that $A$ contains no zeros of $m$. The uniform convergence of $m_n$
to $m$ on $A$ then guarantees that for all sufficiently large $n$, $m_n$
is bounded away from zero in $A$, so $1/m_n$ also converges uniformly
to $1/m$ in $A$. Moreover, Cauchy's integral representation of the derivative
implies that $m_n'$ converges uniformly to $m'$ on $A$. Hence we have
\[ \frac{1}{2\pi i}\int_\Gamma \frac{m_n'(\lambda)}{m_n(\lambda)}d\lambda
 \rightarrow \frac{1}{2\pi i}\int_\Gamma \frac{m'(\lambda)}{m(\lambda)}d\lambda
 \;\; (n\rightarrow\infty),
\]
and the fact that both sides of this equation are integers gives
\be \frac{1}{2\pi i}\int_\Gamma \frac{m_n'(\lambda)}{m_n(\lambda)}d\lambda
 = \frac{1}{2\pi i}\int_\Gamma
\frac{m'(\lambda)}{m(\lambda)}d\lambda \label{eq:4.21} \ee
for all sufficiently large $n$.

From the argument principle, if we let $\nu$ denote the algebraic multiplicity
of $\mu$ as an eigenvalue  -- and hence as pole of $m$ -- we have 
\be -\nu = \frac{1}{2\pi i}\int_\Gamma \frac{m'(\lambda)}{m(\lambda)}d\lambda. 
\label{eq:4.22}\ee
Combining (\ref{eq:4.21}) and (\ref{eq:4.22}) shows that $M_n$ also has
eigenvalues of total algebraic multiplicity $\nu$ inside the contour $\Gamma$, 
for all sufficiently large $n$.
\hfill $\Box$

In section \ref{section:5} we shall give some examples in which the Eastham-Levinson 
asymptotics \cite{kn:Eastham} allow us to verify the hypothesis (\ref{eq:4.19}) and 
hence apply the test for spectral exactness given in Theorems \ref{theorem:4.1}
and \ref{theorem:4.4}.

\section{Pseudospectral inclusion and spectrum-of-sequence inclusion}
\label{subsection:4.4}
In this section we consider a sequence $(T_n)$ of operators on a Hilbert space
$H$. Let $T$ be some other operator on $H$. We denote by $D(T_n)$ the 
domain of $T_n$ and by $D(T)$ the domain of $T$. We shall be interested in
two different types of convergence of $T_n$ to $T$: strong convergence, in
which $T_nf\rightarrow Tf$ for each fixed $f$, and norm resolvent convergence,
in which $\| (\lambda I - T_n)^{-1} - (\lambda I - T)^{-1}\|\rightarrow 0$.
Strong convergence is usually observed when a problem of Sims Case I is
regularized by a sequence of interval truncations, while the stronger
property of norm resolvent convergence is observed when the problem being
regularized is of Sims Case II or III. Strong convergence generally results in a 
very weak type of spectral approximation which is given by Theorem
\ref{theorem:4.5} below; in practice, for differential equation eigenvalue problems,
this is probably not as useful a result as Theorem \ref{theorem:4.4}. 
Norm resolvent convergence, on the other hand, gives a spectral exactness
result (Theorem \ref{theorem:4.7}) which Theorems \ref{theorem:4.1} and 
\ref{theorem:4.2} do not give: Theorems \ref{theorem:4.1} and \ref{theorem:4.2}
do not preclude spectral inexactness in Sims Cases II and III, they merely
give a test for spectral inexactness, whereas Theorem \ref{theorem:4.7} 
precludes spectral inexactness.

The following definition is standard.
\begin{definition}\label{definition:core}
A set ${\cal C} \subseteq D(T)$ is called a core of $T$ if, for every $x\in
D(T)$ and $\epsilon>0$, there exists $x_\epsilon\in {\cal C}$ such that
\[ \| x - x_\epsilon \| < \epsilon, 
   \;\;\; \| Tx - Tx_\epsilon \| < \epsilon. \]
\end{definition}
The following definition is also required.
\begin{definition}
The spectrum of the sequence $(T_n)$, denoted $\sigma(\{ T_n \})$, is the set
\[ \sigma(\{ T_n \}) = \{ \lambda \in \Cc \; | \; \lim_{n\rightarrow\infty} 
\| (\lambda I - T_n)^{-1} \| = +\infty\}. \]
\end{definition}
Note that for selfadjoint operators, $\sigma(\{ T_n \})$ can contain only
points which are limit points of sequences of the form $(\lambda^{(n)})$
in which $\lambda^{(n)}$ lies in the spectrum $\sigma(T_n)$ of 
$T_n$:
\be \sigma(\{ T_n \}) \subseteq \cap_{m\in\Nn}\overline{\cup_{n\geq
m}\sigma(T_n)}. \label{eq:4.23} \ee
To see this, let $\lambda\in \sigma(\{ T_n \})$ and let 
$a_n := \| (\lambda I - T_n)^{-1} \|^{-1}$, so that $a_n\rightarrow 0$
as $n\rightarrow\infty$. By the spectral calculus for the selfadjoint
operator $T_n$, given $\epsilon>0$ there certainly exists a point of 
$\sigma(T_n)$ within distance $a_n+\epsilon$ of $\lambda$: in particular, 
choosing $\epsilon=\frac{1}{n}$ we can find $\lambda^{(n)}\in\sigma(T_n)$
such that $|\lambda - \lambda_n| < a_n + \frac{1}{n}$. Hence for any
integer $m$,
\[ \lambda \in \overline{\cup_{n\geq m}\sigma(T_n)}. \]
This proves (\ref{eq:4.23}). For the non-selfadjoint case the spectral
calculus no longer holds (unless the $T_n$ happen to be normal). 
The following result -- a simple modification of a result of Harrabi 
\cite{kn:Harrabi} -- thus provides for non-selfadjoint operators as
close an analogue of the spectral inclusion result of Reed and Simon 
\cite[Theorem VIII.24]{kn:RS} as is possible, in general.

\begin{theorem}\label{theorem:4.5} 
Let ${\cal C}$ be a core of $T$ and suppose that 
\[ {\cal C} \subseteq \cap_{m\in\Nn}\cup_{n\geq m}D(T_n), \]
so that if $f\in {\cal C}$ then $T_n f $ is defined for all
sufficiently large $n$. 
Let 
\[ \sigma_a(T) := \{ \lambda \in \Cc \; | \; \exists (x_n)_{n\in\Nn}\subseteq
 H \; \mbox{such that}\; \| x_n \| = 1 \; \forall n \; \mbox{and} \;
 \lim_{n\rightarrow\infty} \| (\lambda I - T)x_n \| = 0 \}. \]
Suppose that for each $f\in {\cal C}$ we have $\lim_{n\rightarrow\infty}
\| T_n f - Tf \| = 0$. Then 
\be \sigma_a(T)\subseteq\sigma(\{ T_n \}). \label{eq:4.24} \ee
\end{theorem}
\noindent {\bf Proof}\, Suppose that $\lambda$ does not lie in $\sigma(\{ T_n\})$.
Then there exists $M\in \Rr^{+}$ and a monotone increasing sequence 
$(n_j)_{j\in\Nn}$ such that
\[ \| (\lambda I - T_{n_j})^{-1} \| \leq M. \]
Now let $f\in {\cal C}$. Then $f = (\lambda I - T_{n_j})^{-1}(\lambda I 
 - T_{n_j})f$, whence
\[ \| f \| \leq M \| (\lambda I - T_{n_j})f \| \rightarrow 
 M \| (\lambda I - T)f \| \;\; \mbox{as $j\rightarrow\infty$}. \]
Since this holds for all $f\in{\cal C}$, it follows from Definition
\ref{definition:core} that
\[ \| f \| \leq M \| (\lambda I - T)f \| \;\;\; \forall f\in D(T). \]
Hence $\lambda$ does not lie in $\sigma_a(T)$, which proves the result.
\hfill $\Box$

Because the result (\ref{eq:4.23}) does not hold, in general, for
non-selfadjoint operators, it is not generally possible to
replace $\sigma(\{ T_n\})$ in (\ref{eq:4.24}) by 
$\cap_{m\in\Nn}\overline{\cup_{n\geq m}\sigma(T_n)}$. However, the following result
concerning pseudospectra can be proved.

\begin{theorem}\label{theorem:4.6}
Let $\epsilon>\delta>0$. Let
\be \sigma_\epsilon(T_n) := \{ \lambda \in \Cc \; | \; 
 \| (\lambda I - T_n)^{-1} \| \geq \epsilon^{-1}\}, \label{eq:4.25} \ee
\be \sigma_\delta(T) := \{ \lambda \in \Cc \; | \;
 \| (\lambda I - T)^{-1} \| \geq \delta^{-1}\}. \label{eq:4.26} \ee
Let ${\cal C}$ be a core of $T$ satisfying the same hypotheses as in
Theorem \ref{theorem:4.5}, and suppose that for all $f\in {\cal C}$
we have $\lim_{n\rightarrow\infty}\| T_nf -Tf \| = 0$. Then
\be \sigma_\delta(T) \subseteq \liminf_{n\rightarrow\infty}\sigma_\epsilon(T_n),
\label{eq:4.27} 
\ee
where  
\[ \liminf_{n\rightarrow\infty}\sigma_\epsilon(T_n)
 := \cup_{m\in\Nn}\cap_{n\geq m}\sigma_\epsilon(T_n). \]
\end{theorem}
\noindent {\bf Proof}\, We start by defining 
\be \sigma_\delta(\{ T_n \}) := \{ \lambda \in \Cc \; | \;
\liminf_{n\rightarrow\infty}\| (\lambda I - T_n)^{-1} \| \geq \delta^{-1}\}. 
\label{eq:4.31} \ee
Suppose that $\lambda$ does not lie in $\sigma_\delta(\{ T_n \})$. Then there
exists $\gamma\in (0,\delta^{-1})$ and a monotone increasing sequence 
$(n_j)_{j\in\Nn}$ such that 
\[ \| (\lambda I - T_{n_j})^{-1} \| \leq \gamma < \delta^{-1}. \]
Let $f\in {\cal C}$, and write $f = (\lambda I - T_{n_j})^{-1}
(\lambda I - T_{n_j})f$, giving 
\[ \| f \| \leq \gamma \| (\lambda I - T_{n_j})f\|. \]
Letting $j\rightarrow \infty$ gives
\be  \| f \| \leq \gamma \| (\lambda I - T)f\| \label{eq:4.28} \ee
Eqn. (\ref{eq:4.28}) holds for all $f\in {\cal C}$ and hence, since
${\cal C}$ is a core of $T$, for all $f\in D(T)$. In particular
this implies that $\lambda I - T$ is invertible and
\[ \| (\lambda I - T)^{-1} f \| \leq \gamma \| f \| \]
for all $f\in D(T)$. This implies that $\| (\lambda I - T)^{-1}  \| \leq
 \gamma < \delta^{-1}$, and so $\lambda$ does not lie in
$\sigma_\delta(T)$. We have thus proved
\be \sigma_\delta(T)\subseteq \sigma_\delta(\{ T_n \}). \label{eq:4.29}\ee
The result will be proved if we can show that for $\delta < \epsilon$,
\be \sigma_\delta(\{ T_n \}) \subseteq
\liminf_{n\rightarrow\infty}\sigma_\epsilon(T_n).
\label{eq:4.30} \ee
To do this, suppose that $\mu$ does not lie in 
$\liminf_{n\rightarrow\infty}\sigma_\epsilon(T_n)$. By definition,
\[ \liminf_{n\rightarrow\infty}\sigma_\epsilon(T_n)
 = \cup_{m\in\Nn}\cap_{n\geq m}\{ \lambda \in \Cc \; | \;
 \| (\lambda I - T_n)^{-1} \| \geq \epsilon^{-1} \}, \]
and so for each $m\in \Nn$, $\mu$ does not lie in
\[ \cap_{n\geq m}\{ \lambda \in \Cc \; | \;
 \| (\lambda I - T_n)^{-1} \| \geq \epsilon^{-1} \}. \]
In other words, there exists a subsequence $(T_{n_j})_{j\in\Nn}$
such that 
\[ \| (\mu I - T_{n_j})^{-1} \| < \epsilon^{-1}, \;\;\; 
 j \in \Nn. \]
Hence by definition,
\[ \liminf_{n\rightarrow\infty} \| ((\mu I - T_{n})^{-1} \|
 \leq \epsilon^{-1} < \delta^{-1}. \]
From (\ref{eq:4.31}) we have clearly proved $\mu$ does not lie
in $\sigma_\delta(\{ T_n \}) $. This establishes (\ref{eq:4.30}),
and our proof is complete. \hfill $\Box$

\begin{theorem}\label{theorem:4.7} Let $z\in \Cc$ be fixed.
Suppose that $\| (z I - T_n)^{-1} - (z I - T)^{-1}\|
\rightarrow 0$ as $n\rightarrow\infty$. Then
\be \limsup_{n\rightarrow\infty}\sigma(T_n) \subseteq \sigma(T),
\label{eq:4.32} \ee
where $\limsup$ is defined by
\[ \limsup_{n\rightarrow\infty}\sigma(T_n)
 = \cap_{m\in\Nn}\overline{\cup_{n\geq m}\sigma(T_n)}. \]
\end{theorem}
\noindent {\bf Proof}\, Let $R_n = (zI-T_n)^{-1}$ and let
$R = (zI-T)^{-1}$. Then $\| R_n - R\| \rightarrow 0$ as $n\rightarrow
\infty$. From the results in Kato \cite[IV, \S 3, p. 208]{kn:Kato}
it follows that
\be \limsup_{n\rightarrow\infty}\sigma(R_n) \subseteq \sigma(R). 
\label{eq:4.33} \ee
However the spectrum of $R$ is related to the spectrum of $T$ by
\[ \lambda \in \sigma(T) \;\; \mbox{if and only if} \;\;
   (z-\lambda)^{-1} \in \sigma(R). \] 
[N.B. for our applications, $\infty$ will be an accummulation point
of $\sigma(T)$ and so $0$ will lie in $\sigma(R)$.] A similar relationship
holds between $\sigma(T_n)$ and $\sigma(R_n)$. Thus (\ref{eq:4.33})
implies (\ref{eq:4.32}).
\hfill $\Box$

We shall now show that in Sims Cases II and III, the hypotheses of Theorem 
\ref{theorem:4.7} are satisfied by taking $T_n = M_n$ and $T=M$, where $M_n$ 
and $M$ are the operators of Theorem \ref{theorem:4.1}. Thus Theorem \ref{theorem:4.7}
will supercede Theorem \ref{theorem:4.1} in Sims Cases II and III as a 
guarantee that spectral inexactness is impossible provided the
boundary conditions are correct. Spectral inclusion still holds by
Theorem \ref{theorem:4.4}. Combining all these results will show that
in Sims Cases II and III, provided we generate the $M_n$ using the 
boundary conditions described for these cases in Lemma \ref{lemma:extra}, we
have spectral inclusion and spectral exactness (Theorem \ref{theorem:4.9}
below).

\begin{theorem}\label{theorem:4.8} 
Consider a differential expression of Sims Case II or Case III type at $x=b$.
Let $M$ be a realization of this expression through analytic continuation of
an $m$-function which is a limit of functions $m_n(\cdot)$ for realizations
$M_n$ of the differential operator defined on intervals $[a,b_n]$,
$b_n\nearrow b$ as $n\nearrow\infty$, as described in Lemma \ref{lemma:extra}.
Let $\psi_n(x,\lambda)$ be the 
solution of the differential equation defined by
\[ \psi_n(x,\lambda) = \theta(x,\lambda) + m_n(\lambda)\phi(x,\lambda), \]
and let the $G_n(x,y,\lambda)$ be the Green's functions given by
\[ G_n(x,y,\lambda) = \left\{ \begin{array}{ll}
 -\phi(x,\lambda)\psi_n(y,\lambda), & a < x < y < b, \\
 -\psi_n(x,\lambda)\phi(y,\lambda), & a < y < x < b. \end{array}\right. \]
Let $R_n(\lambda)$ be the extension to $L^2_w[a,b)$ of $(\lambda I-M_n)^{-1}$ defined by
\[ (R_n(\lambda)f)(x) = \int_a^{b_n} G_n(x,y,\lambda)f(y)w(y)dy, 
\;\;\; f\in L^2_w[a,b). 
\]
(see \cite[eqn. (4.2)]{kn:BEMP}). Fix $z\in\Cc$ and suppose that $z$ does not lie 
in the spectrum of $M$ or of any of the $M_n$ for sufficiently large $n$. 
Let $R = (z I-M)^{-1}$. Then $\| R(z) - R_n(z) \| 
\rightarrow 0$ as $n\rightarrow \infty$.
\end{theorem}
\noindent {\bf Proof}\, From \cite[eqn. (4.2)]{kn:BEMP} we know that
\[ (R(z)f)(x) = \int_a^b G(x,y,z)f(y)w(y)dy, \]
in which $\psi(x,z) = \theta(x,z) + m(z)\phi(x,z)$. 
Thus
\[ (R(z)-R_n(z))f)(x) 
 = \int_{b_n}^b G(x,y,z)f(y)w(y)dy
 + \int_a^{b_n}(G(x,y,z)-G_n(x,y,z))f(y)w(y)dy. \]
Because the differential equation is of Sims Case II or Case III, we know that
both $\theta(\cdot,z)$ and $\phi(\cdot,z)$ lie in 
$L^2_w[a,b)$. This implies that
\[ \int_a^b  \hspace{-2mm}w(x)dx \int_a^b w(y)G(x,y,z)^2dy < +\infty, 
 \;\;\; \int_a^b  \hspace{-2mm}w(x)dx \int_a^b w(y)G_n(x,y,z)^2dy < +\infty.
\]
In particular, therefore,
\[ \lim_{n\rightarrow\infty} \int_a^b  \hspace{-2mm}w(x) dx\int_{b_n}^\infty
 |G(x,y,z)|^2 w(y)dy = 0. \]
The bound
\begin{eqnarray*} \| R(z)-R_n(z) \|^2  & \leq & 
   2 \int_a^b \hspace{-2mm}w(x) dx\int_{b_n}^b  \hspace{-2mm}|G(x,y,z)|^2 w(y)dy \\
 &  & + 2 \int_a^b \hspace{-2mm} w(x)dx\int_a^{b}
\hspace{-2mm}w(y)dy|G(x,y,z)-G_n(x,y,z)|^2 \end{eqnarray*}
now yields
\[ \lim_{n\rightarrow\infty} \| R(z)-R_n(z) \|^2  \leq
 2 \lim_{n\rightarrow\infty} 
 \int_a^b  \hspace{-2mm}w(x)dx\int_a^{b}
\hspace{-2mm}w(y)dy|G(x,y,z)-G_n(x,y,z)|^2.\]
We now use the formulae
\[ G(x,y,z) - G_n(x,y,z) = 
\left\{ \begin{array}{ll}
 -(m(z)-m_n(z))\phi(x,z)\phi(y,z), & a < x < y < b, \\
 -(m(z)-m_n(z))\theta(y,z)\phi(y,z), & a < y < x < b, \end{array}\right. \]
to obtain
\[
\lim_{n\rightarrow\infty} \| R(z)-R_n(z) \|^2  \leq
 2 \left(\lim_{n\rightarrow\infty}|m_n(z)-m(z)|^2\right)
 \left(\int_a^b  \hspace{-2mm}w(\xi)|\theta(\xi,z)|^2 d\xi \right)
 \left(\int_a^b  \hspace{-2mm}w(\xi)|\phi(\xi,z)|^2 d\xi \right). 
\]
Since $\lim_{n\rightarrow\infty}|m_n(z)-m(z)|=0$, this proves the
result. \hfill $\Box$
\begin{theorem}(Spectral Inclusion and Spectral Exactness for Sims Cases II
and III). \label{theorem:4.9} \, 
Let $M_n$ and $M$ be as in Theorem \ref{theorem:4.8}. Then
\begin{description}
\item{\em (a)} for every $\lambda$ in the spectrum of $M$, there exists a convergent
sequence $(\lambda^{(n)})_{n\in\Nn}$, with $\lambda^{(n)}$ in the spectrum
of $M_n$, whose limit is $\lambda$;
\item{\em (b)} if $(\lambda^{(n)})_{n\in\Nn}$ is a convergent sequence with
limit $\lambda$ and $\lambda^{(n)}$ lies in the spectrum of $M_n$ for each
$n$, then $\lambda$ lies in the spectrum of $M$.
\end{description}
\end{theorem}
\noindent {\bf Proof} \, By Theorem \ref{theorem:4.8}, the hypotheses of
Theorem \ref{theorem:4.7} are satisfied. This immediately gives (b). 
Turning to (a), we observe that the hypotheses of Theorem
\ref{theorem:4.2}  are satisfied. The result of Theorem \ref{theorem:4.2} 
allows us to use Theorem \ref{theorem:4.3}, which in turn allows us to
use Theorem \ref{theorem:4.4}. The conclusion of Theorem \ref{theorem:4.4}
is precisely (a).\hfill $\Box$

\section{Examples} \label{section:5}
We illustrate the results of the preceding sections with some numerical
examples.

\begin{description}
\item[Example 1]\, An equation of the form
\[ -y'' + c^2 y = \lambda w(x) y, \;\;\; x\in [0,\infty), \]
in which $\Re(c) \neq 0$ and $w(x) = \exp(-3|\Re(c)|x)$, is easily
checked to be in Sims Case II at infinity. Letting $v(x) = 
\exp(-|\Re(c)|x)$ we can define an operator $M$ by 
$(My)(x) = w(x)^{-1}\{ -y'' + c^2 y\}$ for $y\in D(M)$,
where the boundary conditions for $D(M)$ are
\[ y(0) = 0, \;\;\; [y,v](\infty) = 0. \]
The corresponding operator $L$ has domain $D(L)$ specified
by the boundary conditions
\[ y'(0) = 0, \;\;\; [y,v](\infty) = 0. \]
For the operators $M_n$ and $L_n$ on finite intervals
$[0,b_n]$ the boundary conditions at the origin will be
the same as for $M$ and $L$ respectively, while the
boundary condition (\ref{eq:tbc}) at $x=b_n$ will be given,
according to Lemma \ref{lemma:extra}, by
\[ [y,v](b_n) = 0. \]
Using the code described in \cite{kn:SLNSA} we computed the
eigenvalues of the operators $M_n$ and $L_n$ in the box
with corners $100$, $100(1+i)$, $100i$, $0$.  The results,
shown in Table 1, indicate that the eigenvalues of the
$M_n$ and of the $L_n$ converge to distinct points in
the box, and so our test for spectral exactness suggests
that the operator $M$ has eigenvalues close to
$94.4890+i28.8595$ and $24.21335+i14.11108$, while $L$
has eigenvalues close to $3.1163595+i5.808222$ and
$51.51888+i21.3277$. Of course, this is what we would
expect for a Sims Case II problem, by Theorem \ref{theorem:4.9}.

Note that these eigenvalue problems can be formulated as compact 
perturbations of selfadjoint eigenvalue problems, although it is
not immediately clear how one might use this to obtain
spectral inclusion and/or exactness results.

\begin{table}[tbh]
\begin{centering}
\begin{tabular}{|c|c|c|c|c|} \hline 
    &       &                      &                      \\
$n$ & $b_n$ & Eigenvalues of $M_n$ & Eigenvalues of $L_n$ \\
    &       &                      &                      \\ \hline
1   &  5    & 24.21311+i14.10915   & 3.1163619+i5.808222 \\
    &       & 94.4880+i28.8342     & 51.51912+i21.3277 \\ \hline
2   & 10    & 24.21334+i14.11103   & 3.1163595+i5.808219 \\
    &       & 94.4891+i28.8584     & 51.51879+i21.3274 \\ \hline
3   & 15    & 24.21333+i14.11106   & 3.1163591+i5.808221 \\
    &       & 94.4888+i28.8590     & 51.51875+i21.3276 \\ \hline
4   & 20    & 24.21335+i14.11108   & 3.1163595+i5.808222 \\
    &       & 94.4890+i28.8595     & 51.51888+i21.3277 \\ \hline
\end{tabular}
\caption{Example 1 on intervals $[0,b_n]$, using code of 
\protect{\cite{kn:SLNSA}} with $TOL=10^{-7}$.}
\label{table:1}
\end{centering}
\end{table}  

\item[Example 2]\, We consider the (now rather infamous) rotated
 harmonic oscillator problem 
\[ -y'' + c^2 x^2 y = \lambda y, \;\;\; x\in [0,\infty), \;\;\;
 y(0) = 0, \;\;\;
 c\in\Cc, \; \Re(c) > 0,\]
(see Davies \cite{kn:AAD}). This problem
is of Sims Case I at infinity and its eigenvalues are given by
\[ \lambda_k = c(4k+3) \;\;\; k = 0,1,2,\cdots. \]
It is known that the higher index eigenvalues are very
ill-conditioned (when $c^2$ is not positive). Denoting by $M$ the 
operator associated with this problem, this ill-conditioning may 
be explained by the fact that $\| (M - zI)^{-1} \|$ is extremely 
large in very large neighbourhoods of these eigenvalues, making it 
numerically difficult to determine the precise location of the 
poles of $\| (M - zI)^{-1} \|$, which are the eigenvalues.

On the other hand, it is easy to verify that the hypotheses of 
Theorem \ref{theorem:4.3} are satisfied for this problem, so
Theorem \ref{theorem:4.1} still gives a valid test for spectral
inexactness. 

One might expect that this would be of rather academic interest,
given that the operators $M_n$ and $L_n$ on the truncated intervals $[0,b_n]$
will themselves have very ill-conditioned higher index eigenvalues.
To some extent this is correct. However, in Table \ref{table:2} we show the 
result of truncating the interval to $[0,20]$ and locating all the eigenvalues 
in a rectangle in the complex plane with bottom right-hand corner $\lambda=100$ and 
top left-hand corner $\lambda=90i$. The boundary conditions used 
were $y(0) = 0 = y(20)$ for the $M_n$ problem and $y'(0) = 0 = y(20)$
for the $L_n$ problem. The spurious eigenvalues are marked with
asterisks. One can see quite clearly that these eigenvalues 
distinguish themselves by being almost invariant under the
change of boundary condition at the origin. Indeed, in all but
one case the relative differences are less than the tolerance
which was used in the computations ($10^{-5}$). There is also
a problem with `missing' eigenvalues in this table: $M_n$ ought to
have an eigenvalue close to $73 + i 50$ and $L_n$ ought to
have an eigenvalue close to $77 + i 55$, both of which are
missing. Thus the ill-conditioning of these problems may induce
spectral inexactness, for which we seem to be able to test by changing
the boundary conditions, but it can also cause a lack of spectral
inclusion, which is rather more difficult to spot.

\begin{table}[tbh]
\begin{centering}
\begin{tabular}{|c|c|} \hline 
    &                 \\
 Eigenvalues of $M_n$ & Eigenvalues of $L_n$ \\
    &                 \\ \hline
4.3278454 + i 3.1193175  &  1.4426265 + i 1.0397661 \\
10.098296 + i 7.2784044  &  7.2130633 + i 5.1988778 \\
15.868726 + i 11.437476  &  12.983488 + i 9.3579536 \\
21.639117 + i 15.596455  &  18.753885 + i 13.517014 \\
27.409413 + i 19.755403  &  24.524267 + i 17.676070 \\
33.179628 + i 23.914366  &  30.294667 + i 21.835034 \\
38.949823 + i 28.073359  &  36.065062 + i 25.993830 \\
44.720108 + i 32.232282  &  41.835242 + i 30.152542 \\
50.490469 + i 36.391057  &  47.605269 + i 34.311200 \\
56.260702 + i 40.549373  &  53.375120 + i 38.469967 \\
62.032125 + i 44.715888  &  59.146275 + i 42.628491 \\
67.763082 + i 48.620308  &  64.868061 + i 46.792270 \\
70.791811 + i 54.525167* &  70.836011 + i 54.242901* \\
76.994286 + i 49.224379  &  71.821431 + i 49.862180 \\
72.268485 + i 64.384873* &  72.268524 + i 64.384384* \\
73.809759 + i 74.921450* &  73.809759 + i 74.921449* \\
75.474904 + i 86.054406* &  75.474905 + i 86.054405* \\
87.360734 + i 47.089232  &  82.109895 + i 48.176485 \\
98.348465 + i 44.849646  &  92.771537 + i 45.977747 \\ \hline
\end{tabular}
\caption{Testing for spectral inexactness on the rotated
 harmonic oscillator}
\end{centering}
\label{table:2}
\end{table}

\item[Example 3]\, Consider the problem of locating resonances of
the equation
\be -y'' + 16x^2 \exp(-x)y = \lambda y, \;\;\; y'(0) = 0, \;\;\;
 x\in [0,\infty). \label{eq:reseq1} \ee
Using the `complex scaling' method, the resonances of this
problem are of the form $\mbox{e}^{-2i\theta}\mu_\theta$ where $\mu_\theta$
is an eigenvalue of the non-selfadjoint problem
\be -z'' + 16x^2\mbox{e}^{2i\theta} \exp(-x\mbox{e}^{i\theta})z = \mu z, \;\;\;
z'(0) = 0, \;\;\; x\in [0,\infty), \label{eq:reseq2} \ee
(see Hislop and Segal \cite{kn:Hislop}). The rotation angle
$\theta> 0$ must be such that the function $x\mapsto 16x^2\mbox{e}^{2i\theta}
\exp(-x\mbox{e}^{i\theta})$ lies in $L^1[0,\infty)$, and in particular
therefore $\theta < \pi/2$. Resonances have the property that 
$\mbox{e}^{-2i\theta}\mu_\theta$ is independent of $\theta$, so in general
not all eigenvalues of (\ref{eq:reseq2}) yield resonances: one should
carry out the computations for at least two different values of $\theta$
to identify resonances.

In addition to the complications caused by the fact that some
eigenvalues of (\ref{eq:reseq2}) do not correspond to resonances,
we have the additional problem that (\ref{eq:reseq2}) is a singular
problem and must be regularized by interval truncation. This truncation
process might introduce spurious eigenvalues, not corresponding to 
eigenvalues of (\ref{eq:reseq2}), which we need to be able to detect. 
Theorem \ref{theorem:4.1} gives a way to do this.

Using a rotation angle $\theta=1.1$ and a truncated interval
$[0,100]$ with boundary condition $y(100)=0$ we computed both
the eigenvalues of the equation in (\ref{eq:reseq2}) with
$y'(0) = 0$ and the eigenvalues for the same equation but with
the boundary condition $y(0) = 0$ using the code of \cite{kn:SLNSA}.
We asked the code to find, in the $\mu$ plane, all the eigenvalues
in the box with corners $(-0.01,0.01)$, $(-0.01,5)$, $(-10,5)$,
$(-10,0.01)$, with a tolerance of $10^{-6}$. The results, rotated 
back into the $\lambda$ plane via $\lambda = \mbox{e}^{-2i\theta}\mu$, 
are shown in Table \ref{table:3}. 

\begin{table}[tbh]
\begin{centering}
\begin{tabular}{|c|c|} \hline 
 & \\
    Alleged resonances & Alleged resonances                \\
      with $y(0)=0=y(100)$ & with $y'(0)=0=y(100)$ \\
    on $[0,100]$ &              on $[0,100]$   \\ 
 & \\ \hline
2.429823932+i 2.95502902 & 2.429823937 +i 2.95502903 \\
3.869964809-i 0.74439879 & 3.869964804 -i 0.74439879 \\
0.554661821+i 0.66915540 & 0.554661961 +i 0.66915556 \\
                         & 2.861786706 -i $1.6 \times 10^{-6}$ 
 \\ \hline
\end{tabular}
\caption{Testing for spurious resonances due to interval truncation}
\end{centering}
\label{table:3}
\end{table}

Of the four alleged resonances found with $y'(0)=0$, three
are virtually unchanged when the boundary condition is changed
to $y(0)=0$. Theorem \ref{theorem:4.1} indicates that these
are probably spurious. This is obvious for the two which have
positive imaginary parts, as resonances lie in the lower
half plane by definition; however, without Theorem \ref{theorem:4.1}
it would not have been obvious for the alleged resonance at 
$3.8699648-i0.7443988$. For this problem we believe that the only 
genuine resonance found, for boundary condition $y'(0)=0$, is 
the one at $2.861786706 -i 1.6 \times 10^{-6}$. In fact, of the
four alleged resonances this is the only one which is invariant
under a change of the rotation angle $\theta$; however, in general
it is not clear that a spurious resonance generated by interval 
truncation would always fail to be invariant under change of
$\theta$.

\end{description}

\bibliographystyle{plain}

\end{document}